\documentclass{article}
\usepackage[utf8]{inputenc}

\usepackage[a4paper, margin=3cm]{geometry}
\usepackage{amsmath}
\usepackage{amssymb}
\usepackage{bbm}
\usepackage{amsthm}
\usepackage{graphicx}
\usepackage[title]{appendix}

\newtheorem{tw}{Theorem}
\newtheorem{lem}[tw]{Lemma}
\newtheorem{prop}[tw]{Proposition}
\newtheorem{de}[tw]{Definition}
\newtheorem{cor}[tw]{Corollary}

\theoremstyle{remark}

\newcommand{\dd}{\mathrm{d}}
\newcommand{\ddiv}{\mathrm{div}\,}

\title{Weak solutions for the Stokes system for compressible fluids with general pressure}
\author{Maja Szlenk}
\date{\today}

\numberwithin{equation}{section}
\numberwithin{tw}{section}

\begin{document}

\maketitle

\begin{abstract}
    We prove existence and uniqueness of global in time weak solutions for the Stokes system for compressible fluids with a general, non-monotone pressure. We construct the solution at the level of Lagrangian formulation and then define the transformation to the original Eulerian coordinates. For nonnegative and bounded initial density the solution is also nonnegative for all $t$ and belongs to $L^\infty([0,\infty)\times\mathbb{T}^d)$. A key point of our considerations is the uniqueness of such transformation. Since the velocity might not be Lipschitz continuous, we develop a method which relies on the results of Crippa \& De Lellis, concerning regular Lagriangian flows. The uniqueness is obtained thanks to the application of a certain weighted flow and detail analysis based on the properties of the $BMO$ space.
\end{abstract}


\section{Introduction}

The Stokes system is an approximation of the Navier-Stokes equations for small Reynolds number. In such cases the advective intertial forces are relatively small and explicit dependence on time and convective term can be omitted. This is a typical situation for highly viscous fluids or when the velocities are very small. The flows satisfying these conditions are called Stokes or creeping, and they occur in numerous biological and physical problems, e.g. to describe dynamics of the blood in a process of sedimentation \cite{blood}, or to model swimming of microorganisms \cite{swimming2,swimming,microorganisms}. Other applications include also engineering, where the Stokes flow is used in the process of designing microfluids and microdevices \cite{microfluids2,microfluids}. The Stokes model is also connected to the Darcy law, which describes the flow of a fluid through porous media. Such phenomena are observed in biological tissues \cite{tumor1,tumor2} and have many applications in petroleum engineering \cite{petroleum, petroleum2}. The other situation, where the fluid motion is governed by the Stokes equation is a laminar flow. In this case the fluid particles move in adjacent layers, with little mixing between them. 

In this article we consider the compressible Stokes flow on the $d$-dimensional torus $\mathbb{T}^d$
\begin{equation}\label{NS}
\left\{ \begin{aligned} &\varrho_t + \text{div}(\varrho u) = 0,  \\ 
 &-\mu\Delta u -\nabla(\lambda+\mu)\ddiv u+\nabla p(\varrho) = 0, \end{aligned}\right.  \end{equation}
 where $\varrho\colon[0,T]\times\mathbb{T}^d\to\mathbb{R}$ and $u\colon[0,T]\times\mathbb{T}^d\to\mathbb{R}^d$ are the sought fluid density and velocity field. The function $p(\varrho)$ denotes the pressure term, and the parameters $\mu,\lambda$ represent the first and the second viscosity. \\
 
 We further assume that the flow is irrotational, namely $\mathrm{rot}u=0$. It is equivalent to velocity having the structure of a gradient flow, namely $u(t,x)=\nabla\phi(t,x)$ for some $\phi\colon [0,T]\times\mathbb{T}^d\to\mathbb{R}$. In consequence, we obtain the condition $\displaystyle \int_{\mathbb{T}^d} u(t,x)\dd x=0$. The second equation of (\ref{NS}) can be then rewritten in terms of effective viscous flux, which turns out to be constant. Therefore instead of the second equation of (\ref{NS}) we obtain
 \[ (\lambda+2\mu)\ddiv u=p(\varrho)-\{p(\varrho)\}, \]
 where $\displaystyle\{f\}=\frac{1}{\mathbb{T}^d}\int_{\mathbb{T}^d}f(x)\dd x$. As the qualitative properties of solutions do not depend of the values of $\lambda$ and $\mu$, without loss of generality we take $\lambda+2\mu=1$. Under these assumptions the system (\ref{NS}) can be transformed into
 \begin{equation}\label{NS2}
 \begin{aligned}
 \varrho_t &+ \text{div}(\varrho u) = 0, \\
 \ddiv u &=p(\varrho)-\{p(\varrho)\}.
 \end{aligned}
 \end{equation}
 
 Equation (\ref{NS2}) is coupled with the initial condition on the density, which is assumed to be bounded and nonnegative, namely
 \[ \varrho_{|_{t=0}} := \varrho_0 \in L^\infty(\mathbb{T}^d), \quad \varrho_0\geq 0. \]
 In particular, we do not require the density to be strictly positive, hence for example $\varrho_0$ can be equal $\mathbbm{1}_A$ for some $A\subset\mathbb{T}^d$. \\
 
 Our method allows the pressure to be in a quite general form. We require $p(\varrho)$ to be of class $C^1$ and unbounded, so that in particular we can choose a sequence $\varrho_n\to\infty$ such that $p(\varrho_n)\to\infty$. Moreover, we assume that there exist constants $C, C_1,C_2,\bar{\varrho}$ such that $p$ satisfies the inequality
 \begin{equation}\label{P} -C \leq p(\varrho) \leq CP(\varrho) := C\left(\varrho\int_{\bar{\varrho}}^\varrho \frac{p(s)}{s^2} \, \dd s +C_1\varrho+C_2\right). \end{equation} 
 Observe that as $p(\varrho)$ is bounded from below, we can choose $\bar{\varrho},C_1,C_2$ such that $P(\varrho)\geq 0$ for all $\varrho\geq 0$, hence from now on we will assume that $P(\varrho)\geq 0$. The properties of functions satisfying (\ref{P}) are discussed in subsection \ref{press_sect}.

\subsection{Statement of the main theorem.}
The mathematical theory of weak solutions to compressible fluid equations has been widely developing in the last twenty years, since the groundbreaking results of Lions in 1998 \cite{lions} and Feireisl \cite{feireisl,feireisl2} in 2001. They proved the existence of weak solutions to compressible Navier-Stokes equations, provided that the pressure term is of the form $p(\varrho)=\varrho^\gamma$ with $\gamma>\frac{9}{5}$ and $\gamma>\frac{3}{2}$ respectively. In \cite{non-monotone} and \cite{non-monotone2}, this method was also adjusted to the pressure which is non-monotone on some finite interval. In particular, it allows to deal with the pressures expressed via equations of state, which are of more complex form than ideal gas, the model example being the van der Waals' equation of state. The Lions \& Feireisl approach can be also adjusted to more complex systems, for example Navier-Stokes-Fourier system \cite{feireisl-fourier, feireisl-mucha} and other including entropy transport \cite{mucha-entropy} or heat conductivity \cite{bresch-desjardins}. 

 The important results concerning non-monotone pressure laws are due to D. Bresch and P.-E. Jabin \cite{BJ2,bresch-jabin}. Their method, based on the Kolmogorov compactness criterion, allows to deal with the pressure satisfying 
  \begin{equation}\label{BJ} C^{-1}\varrho^\gamma-C \leq p(\varrho) \leq C\varrho^\gamma+C 
\quad \text{ and } \quad |p'(\varrho)|\leq \bar{p}\varrho^{\gamma-1} \end{equation}
for $\gamma>\frac{9}{5}$. In context of our work, we refer the reader especially to \cite{bresch-jabin}, where they presented their approach also on a modification of the Stokes system. In this case they proved the existence of global weak solutions for any $\gamma>1$ with the same regularity as in the isentropic case, namely $\varrho\in L^\infty(0,T;L^\gamma)$. 

Besides the relaxed conditions on the pressure, the Bresch \& Jabin compactness criterion can be also applied to various classes of equations, where the Lions and Feireisl method was insufficient. One of the examples are the systems with the additional term in the continuity equation. Such models can be obtained from the multi-fluid systems \cite{mucha} or appear in the mathematical modelling of tumor growth \cite{vauchelet-zatorska}. In this case the additional term results in lack of compensated compactness between $\textrm{div}u$ and the pressure and therefore the classical method fails. However, the Bresch \&
Jabin criterion allows to omit this problem. The problem with the convergence of effective viscous flux arise also in the anisotropic case. The recent result \cite{bresch-anisotropic}, concerning anisotropic compressible Stokes system, omits this problem by controlling a certain defect measure associated to the pressure.

The other related topic are steady flows, where the behaviour of the fluid does not depend on time. The equations describing such flows are the classical equations of fluid mechanics with the time derivative set as zero, for example compressible Navier-Stokes \cite{piasecki_steady,plotnikov} or Navier-Stokes-Fourier system \cite{mucha-pokorny, pokorny-novotny}. Again, in the case of the more general pressure laws, the classical Lions-Feireisl method cannot be directly applied to obtain a weak solution. One of the ways to cope with that problem was introducing the variational entropy solutions \cite{entropy1, entropy2}, where the total energy balance is replaced with entropy inequality. Another analyzed system is the steady Oseen flow \cite{piasecki_oseen}, which is a linearization of the Navier-Stokes system with partial consideration of the convective forces. Note that in our case the explicit dependence on time is removed only in the momentum equation, therefore this system can be consider as an intermediate step between steady and unsteady flows. \\
 
 As the Stokes system has a simpler structure than the Navier-Stokes, the analysis can be carried out in the more general setting. If $\varrho_0$ has higher regularity, then the solutions to the Stokes problem exist and are unique for the general $p$ (see \cite{lions}, Remark 8.14). However, in case of $\varrho_0\in L^\infty(\mathbb{T}^d)$, the classical method requires the monotonicity condition on the pressure and the uniqueness was not established. In this paper we obtain the uniqueness of solutions to the Stokes system in case of the low regularity of the initial density and under very general pressure laws.  Our main theorem states as follows:
\begin{tw}\label{main}
 Let $\varrho_0\in L^\infty(\mathbb{T}^d)$, $\varrho_0\geq 0$ and the pressure satisfy (\ref{P}). Then there exists a unique global in time weak solution to (\ref{NS2}), satisfying
 \[\begin{aligned} \varrho &\in L^\infty([0,\infty)\times\mathbb{T}^d), \\
 u &\in L^\infty([0,\infty)\times\mathbb{T}^d), \\ 
 \nabla u &\in L^\infty([0,\infty);BMO), \; \ddiv u\in L^\infty([0,\infty)\times\mathbb{T}^d).
 \end{aligned}\]
\end{tw}

Our approach is based on the Lagrangian reformulation of the system, which allows us to obtain a global $L^\infty$ estimate on the density. Having that estimate, we can straightforwardly apply Bresch and Jabin method to obtain compactness, and in consequence existence of solutions besides the relaxed conditions on the pressure. However, using the results from the theory of transport equations and classical harmonic analysis, we were able to establish also uniqueness of solutions. \\

The first step is to show the existence of a unique solution in the Lagrangian formulation. Then, it suffices to define a transformation back to the original Eulerian coordinates. However, there appear some difficulties in the construction, concerning low regularity of the solution of the Lagrangian reformulation. The obtained regularity of $u$ provides that $\nabla u\in L^\infty(0,T;BMO)$. In particular, $\nabla u$ may not be bounded in $x$ and in consequence the flow $x(t,y)$ generated by $u$ may not be invertible on the whole torus. Nevertheless, the divergence of $u$ remains bounded, therefore $x(t,y)$ is a regular Lagrangian flow and we can treat it using the properties from \cite{colombo-crippa} and \cite{crippa-delellis}. The key tool needed in our analysis is the recent result of Crippa and De Lellis \cite{crippa-delellis} concerning the stability of regular Lagrangian flows in $L^1$, which allows us to pass to the limit with the smooth approximation of the system. 

The above reasoning, however, does not preserve uniqueness. We need to prove the latter by taking two solutions $u_0, u_1$, which coincide in Lagrangian coordinates, and show that $u_0=u_1$. For this purpose we introduce a certain weighted flow between $u_0$ and $u_1$ satisfying
\[
\dot{x}_s = su_1(t,x_s)+(1-s)u_0(t,x_s), \quad x_s(0,y) = y. \]
Having the crucial $BMO$ regularity of $\nabla u_i$, $i=1,2$, we are able to use the John-Nirenberg inequality and a certain integral inequality for functions of bounded mean oscillation, proved in \cite{mucha-rusin, mucha-transport}. Then, carefully combining the estimates we show that $x_s$ does not depend on $s$ and in consequence $u_0=u_1$. \\

The rest of the paper is divided into sections, which contain the main steps of the proof of Theorem \ref{main}. The structure of the proof is as follows:
\begin{itemize}
    \item In section \ref{lagr_sect} we present the a priori estimates and results at the level of Lagrangian coordinates, namely the $L^\infty$ bounds and unique existence of a solution in the Lagrangian reformulation.
    \item Section \ref{uniq_sect} contains the necessary tools and definitions, together with the proof of the uniqueness of the solutions to equation (\ref{NS2}).
    \item In section \ref{euler_sect} we define the transformation from Lagrangian to Eulerian coordinates using the construction from \cite{mucha}, and therefore prove the existence of solutions to (\ref{NS2}). Note that the estimates obtained in section \ref{lagr_sect} provide the $BMO$ regularity of the gradient, necessary to obtain uniqueness, and the existence of solutions can be done independently, using for example standard Lions method and the Bresch \& Jabin compactness criterion. However, we present here an alternative approach, which is set in the framework of Lagrangian regular flows and therefore is more consistent with the rest of this paper.
\end{itemize}

\paragraph{Notation remarks:} For the notational simplicity, we omit the subscript while integrating over the torus, namely
\[ \int \dd x \ := \ \int_{\mathbb{T}^d} \dd x. \]
By $\{\cdot\}_Q$ we denote the mean integral over $Q$, while in the case of the whole torus we again omit the subscript. \\
Moreover, to distinguish between the norms in $L^\infty(\mathbb{T}^d)$ and $L^\infty([0,T]\times\mathbb{T}^d)$, we denote
\[ \|\cdot\|_\infty :=\|\cdot\|_{L^\infty(\mathbb{T}^d)} \ \text{ and } \ \|\cdot\|_{\infty,T} :=\|\cdot\|_{L^\infty([0,T]\times\mathbb{T}^d)}. \]

 \subsection{Discussion on the pressure}\label{press_sect}
 
 Observe that our condition in particular includes the assumptions from \cite{bresch-jabin}:
     if $p$ satisfies (\ref{BJ}), then $P(\varrho)\geq c\varrho^\gamma-c\varrho$ and for a sufficiently large $C$ we obtain (\ref{P}). However, our assumptions also allow the pressure to drop to $0$ for arbitrary large $\varrho$ and we do not require the bounds on the derivative. 
 
 Such class of possible pressures includes many physical situations. It contains the cases which were covered before, for example van der Waals' fluid. Moreover, the admissible pressures also can be expressed via virial expansion, namely defined as a power series of the density:
 \[ p(\varrho)=\sum_{k=1}^\infty B_k\varrho^k, \]
where coefficients $B_k$ depend on the temperature and are derived from statistical mechanics. The virial equation of state was also considered in \cite{BJ2}, but our result allows wider range of pressures of this type. The other case, where our result may be applicable, is the use in biological models, where the pressure term is responsible for interactions between different biological agents and therefore can take form other than resulting from physical constitutive laws.
 
 Let us present some further properties of $p$ satisfying (\ref{P}):
 \begin{itemize}
 \item Condition (\ref{P}) implies that in particular $p(\varrho)\leq \tilde{C}\varrho^\gamma+\tilde{C}$ for $\varrho\geq \bar{\varrho}$ and some $\gamma>1$: \\

     Let $p$ satisfy (\ref{P}) and define $\alpha(\varrho)=\int_{\bar{\varrho}}^\varrho \frac{p(s)}{s^2}\dd s$, $\alpha(\bar{\varrho})=0$. Then 
     \[ \alpha'(\varrho)=\frac{p(\varrho)}{\varrho^2}\leq \frac{C\alpha(\varrho)}{\varrho}+\frac{C}{\varrho}\left(1+\frac{1}{\varrho}\right) \]
     and by the comparison criterion $\alpha(\varrho) \leq \bar{C}\varrho^C -\frac{C}{C+1}\frac{1}{\varrho}-1$ for $\varrho\geq\bar{\varrho}$, where $\bar{C}$ depends on $\bar{\varrho}$. Therefore
     \[ p(\varrho)\leq C(\varrho\alpha(\varrho)+\varrho+1) \leq \bar{C}C\varrho^{C+1}+\frac{C^2}{C+1}\leq \tilde{C}\varrho^\gamma+\tilde{C} \]
     for $\gamma=C+1$ and a suitable $\tilde{C}$.
     \item On any finite interval we can estimate $p$ by sufficiently large constant, hence in particular (\ref{P}) is fulfilled. Therefore to check if indeed $p$ satisfies (\ref{P}) for all $\varrho\geq 0$, it suffices to analyse the behaviour of $P$ when $\varrho\to\infty$. It is immediate to check that if $P$ satisfies 
     \[ \liminf_{\varrho\to\infty} \frac{P(\varrho)}{\varrho^{\gamma}}=\liminf_{\varrho\to\infty}\frac{1}{\varrho^{\gamma-1}}\int_{\bar{\varrho}}^\varrho \frac{p(s)}{s^2}\dd s \geq c >0 \]
     for some $\gamma>1$, then $p$ satisfies (\ref{P}), however these conditions are not equivalent.
     \item  The most significant difference between our class of admissible pressures and the cases considered before is that we allow the pressure to drop to $0$ even for large $\varrho$. Moreover, the derivative of $p$ may grow arbitrarily fast. For example, let $f\colon [0,\infty)\to [0,\infty)$ be a smooth, increasing function such that $f'$ is also increasing and define 
     \[ p(\varrho)=\varrho^2(1+\cos(f(\varrho))). \]
     By the alternating series test, the integral
     \[ \int_{\bar{\varrho}}^\infty \cos(f(x))\dd x = \int_{f(\bar{\varrho})}^\infty \frac{\cos(y)}{f'(f^{-1}(y))}\dd y \]
     is convergent in the sense of Riemann. Therefore we have 
     \[ \lim_{\varrho\to\infty}\frac{1}{\varrho}\int_{\bar{\varrho}}^\varrho \frac{p(s)}{s^2}\dd s = \lim_{\varrho\to\infty}\frac{1}{\varrho}\int_{\bar{\varrho}}^\varrho 1+\cos(f(s))\dd s = 1 \]
     and $p$ satisfies (\ref{P}). Moreover, it periodically drops to $0$ and the derivative of $p$ depends on $f'$, which can be arbitrarily large. 
     \item The condition $p(\varrho)\leq C\varrho^\gamma+C$ is not sufficient to obtain (\ref{P}). For example, let $\eta$ be a smooth function supported in $[-1,1]$ such that $0\leq \eta\leq 1$ and $\eta(0)=1$. Define
     \[ p(\varrho) = \left\{ \begin{aligned}
     &\varrho^2\eta(2^k(\varrho-k)) \; \text{ for } k-2^{-k} \leq \varrho \leq k+2^{-k}, k=1,2,... \\
     &0 \; \text{ otherwise}.
     \end{aligned}\right. \]
     Then 
     \[\begin{aligned} \int_0^k \frac{p(s)}{s^2}\dd s &= \sum_{i=1}^{k-1}\int_{i-2^{-i}}^{i+2^{-i}}\eta(2^i(s-i))\dd s+\int_{k-2^{-k}}^k \eta(2^k(s-k))\dd s \\
     &= \sum_{i=1}^{k-1} 2^{-i}\int_{-1}^1 \eta\dd x + 2^{-k}\int_{-1}^0 \eta\dd x = (1-2^{-k+1})\int_{-1}^1\eta\dd x + 2^{-k}\int_{-1}^0\eta \dd x \end{aligned}\]
     and therefore for any $C$ we can choose sufficiently large $k$ such that 
     \[ p(k)=k^2 \geq C\left(k\int_0^k \frac{p(s)}{s^2}\dd s+k+1\right) \approx \tilde{C}k. \]
     \end{itemize}
     
 
 \section{The $L^\infty$ bound on the density}\label{lagr_sect}
 
  \subsection{Energy estimates}
 First, we obtain the a priori estimates for our solutions.
 \begin{lem}
 Let $(\varrho,u)$ be sufficiently smooth solutions to (\ref{NS2}). Then they satisfy the estimate
 \[ \int_0^T\int (\ddiv u)^2 \; \dd x \dd t + \sup_{t\in[0,T]} \int P(\varrho) \; \dd x \leq C\int P(\varrho_0) \; \dd x. \]
 \end{lem}
 Proof: By multiplying (\ref{NS2})\textsubscript{2} by $\ddiv u$ and integrating on torus, we get
 \[ \int (\ddiv u)^2 \; \dd x - \int p(\varrho)\ddiv u \; \dd x = 0. \]
Moreover,
\[\begin{aligned} -\int p(\varrho)\ddiv u \; \dd x &= \int\nabla p\cdot u \; \dd x = \int \frac{p'(\varrho)}{\varrho}\nabla\varrho\cdot(\varrho u) \; \dd x = -\int P'(\varrho)\cdot \text{div}(\varrho u) \; \dd x \\
&= \int P'(\varrho)\varrho_t \; \dd x = \frac{\dd}{\dd t}\int P(\varrho) \; \dd x. \end{aligned} \]
 After integration over time, we get the desired estimate. \qed \\

The above a priori bounds also provide that $p(\varrho)\in L^\infty(0,T;L^1)$. Indeed, from assumption (\ref{P}) we have
\begin{equation}\label{p_oszac}
 \sup_{t\in [0,T]}\int p(\varrho)\dd x \leq C\sup_{t\in[0,T]}\int P(\varrho)\dd x \leq C. \end{equation}
 
 \subsection{The Lagrangian formulation}
 To prove the global $L^\infty$ estimate on the density, we need to rewrite the equation (\ref{NS2}) in Lagrangian coordinates, that is, we carry out a certain change of variables, which allows us to reduce the continuity equation to a simple ordinary differential equation. 
 
 If $u$ is a velocity field, then the trajectory of a single fluid parcel moves along the integral curve of $\dot{x} = u(t,x)$. Therefore, if at the starting time the particle was at a point $y$, then after the time $t$ it would be at the point $x(t,y)$, where $x(t,y)$ is the solution to the Cauchy problem
 \[\begin{aligned} \dot{x}(t,y) &= u(t,x(t,y)), \\
 x(0,y) &= y. \end{aligned} \]
 Moreover, differentiating this equation with respect to $y$ and using Liouville formula, we obtain the equation for the Jacobian $J$ of $\frac{\dd x}{\dd y}$:
 \[ \begin{aligned} \dot{J}(t,y) &= \ddiv u(t,x(t,y))J(t,y), \\
 J(0,y) &= 1, \end{aligned} \]
 and in consequence $\displaystyle J(t,y)=\exp\left(\int_0^t \ddiv u(s,x(s,y))\dd s\right)$. 
 
 We rewrite the system (\ref{NS2}) in terms of $y$ coordinates instead of $x$, so at a time $t$ the space variable is the position of a parcel starting from $y$. Let
 \[ \eta(t,y)=\varrho(t,x(t,y)) \; \; \; \text{ and } \; \; \; \sigma(t,y)=\ddiv u(t,x(t,y)). \]
 Then the equation (\ref{NS2}) has the form
 \begin{equation}\label{eta}
 \begin{aligned} \partial_t\eta &+\eta\sigma=0, \\
 \sigma &=p(\eta)-\{p(\eta)\}_\sigma, \end{aligned}
 \end{equation}
 where $\{\cdot\}_\sigma$ is the mean integral in the new coordinates given by
 \[ \{f\}_\sigma = \frac{1}{|\mathbb{T}^d|}\int f(t,y)\exp\left(\int_0^t\sigma(s,y)\mathrm{d}s\right)\mathrm{d}y. \]
 
 Now for the system (\ref{eta}) we obtain the following result:
 \begin{tw}\label{lagrange}
  For $\varrho_0\in L^\infty(\mathbb{T}^d)$ and any $T>0$ there exists a unique solution 
  \[ (\eta, \sigma)\in L^\infty([0,T]\times\mathbb{T}^d)\times L^\infty([0,T]\times\mathbb{T}^d) \]
  to the equation (\ref{eta}) with the initial condition $\eta(0,y)=\varrho_0(y)$. Moreover, there exists a constant $r$, independent of $T$, such that 
  \[ \|\eta\|_{\infty,T} \leq r. \]
 \end{tw}
 
 The immediate consequence of Theorem \ref{lagrange} is the similar uniform bound on $\sigma$. As $T$ is arbitrary, we hence obtain the existence of a solution on the whole real line. \\
 
 The proof of the existence of a unique solution is a standard application of the Banach fixed point theorem and is presented in the Appendix \ref{existence_proof}. Here we show the second part of Theorem \ref{lagrange}, namely the $L^\infty$ bound on $\eta$. 
 \begin{prop}
 If $(\eta,\sigma)\in L^\infty([0,T]\times\mathbb{T}^d)\times L^\infty([0,T]\times\mathbb{T}^d)$ is a solution to (\ref{eta}), then
 \[ \|\eta\|_{\infty,T} \leq r, \]
 where the constant $r$ does not depend on $T$.
 \end{prop}
 
Proof:  It is useful to solve the first equation of (\ref{eta}). We get the identity
\begin{equation}\label{wzor_eta} \eta(t,y)=\varrho_0(y)\exp\left(-\int_0^t\sigma(s,y)\mathrm{d}s\right). \end{equation}

Note that this explicit formula for $\eta$ provides that in particular if $\varrho_0$ is strictly positive, then for any $t$ the density is strictly positive as well. \\

Using (\ref{wzor_eta}), we see that $\eta$ is also continuous with respect to time. For a fixed $y$, we have
\[\begin{aligned} |\eta(t+\varepsilon,y)-\eta(t,y)| &= \varrho_0(y)\exp\left(-\int_0^t \sigma(s,y)\dd s\right)\left|\exp\left(-\int_t^{t+\varepsilon}\sigma(s,y)\dd s\right)-1\right| \\
&\leq \varrho_0(y)e^{t\|\sigma\|_{\infty,T}}\left|e^{\varepsilon\|\sigma\|_{\infty,T}}-1\right|. \end{aligned} \]
Hence $|\eta(t+\varepsilon,y)-\eta(t,y)|$ goes to $0$ as $\varepsilon\to 0$ and indeed $\eta(\cdot,y)$ is continuous. \\

 The continuity of $\eta$ allows us to show global boundedness. Recall that by virtue of energy estimates, the mean value $\{p(\eta)\}_\sigma$ is bounded. Let
 \[ \sup_{t\in[0,T]}\{p(\eta)\}_\sigma=M. \]
 From (\ref{p_oszac}) $M$ is finite and does not depend on $T$. As $p$ is unbounded, we can choose $r>\|\varrho_0\|_{L^\infty(\mathbb{T}^d)}$ such that $p(r)>M$. Then at the point $\eta=r$ we get
 \[ \partial_t\eta|_{\eta=r} =-r(p(r)-\{p(\eta)\}_\sigma) < -r(p(r)-M) <0. \]
 However, as for a fixed $y$ $\eta(\cdot,y)$ is continuous and $\eta(0,y)=\varrho_0(y)<r$, if it exceeds the value $r$, it must have a nonnegative derivative at that point, which gives a contradiction. Hence for any $y\in \mathbb{T}^d$ the function $\eta(\cdot,y)$ is also bounded by $r$. \qed
  
 \section{Uniqueness of solutions}\label{uniq_sect}
 Using Lagrangian coordinates introduced in the previous section, we are able to show that the solutions to (\ref{NS2}) are unique.
 
 \begin{tw}
 If $(\varrho_i, u_i), \; i=1,2$ are solutions to (\ref{NS2}) with the regularity from Theorem \ref{main}, then $(\varrho_1,u_1)=(\varrho_2,u_2)$.
 \end{tw}
 
First, let us show that if $(\varrho, u)$, is a solution to (\ref{NS2}) satisfying
\[ \varrho, u, \ddiv u \in L^\infty([0,T]\times\mathbb{T}^d), \quad \nabla u \in L^\infty(0,T;BMO), \]
then in the Lagrangian coordinates it satisfies (\ref{eta}). From the classical results from transport theory of DiPerna \& Lions \cite{diperna-lions} it follows that for such $u$ there exists unique flow $x(t,y)$, such that $x\in C(0,T;L^p)$ for any $1<p<\infty$ and
\[ \dot{x} = u(t,x), \quad x(0,y)=y. \]
Moreover, if $\varrho$ is a solution to the continuity equation
\begin{equation}\label{continuity} \varrho_t+\ddiv(\varrho u) =0, \quad \varrho(0,\cdot)=\varrho_0, \end{equation}
then the function $\varrho(t,x(t,y))$ is equal to
\[ \varrho(t,x(t,y)) = \varrho_0(y)\exp\left(-\int_0^t \text{div}u(s,x(s,y))\dd s\right) = \varrho_0(y)\exp\left(-\int_0^t \sigma(s,y)\dd s\right) \]
and therefore $\eta(t,y)=\varrho(t,x(t,y))$ satisfies the first equation of (\ref{eta}) with $\sigma(t,y)=\ddiv u(t,x(t,y))$. Furthermore, taking the second equation of (\ref{eta}) at a point $x(t,y)$, we obtain
\[ \sigma(t,y) = p(\eta(t,y))-\frac{1}{\mathbb{T}^d}\int p(\varrho(t,x))\dd x. \]
However, by Lemma 3.1. from \cite{colombo-crippa}, for any $f\in L^1(\mathbb{T}^d)$ we have 
\[ \int f(x) \dd x = \int f(x(t,y))e^{\int_0^t \ddiv u(s,x(s,y))\dd s}\dd y = \int f(x(t,y))e^{\int_0^t \sigma(s,y)\dd s}\dd y. \]
Thus 
\[ \int p(\varrho(t,x))\dd x = \int p(\varrho(t,x(t,y))e^{\int_0^t \sigma(s,y)\dd s}\dd y = \int p(\eta(t,y))e^{\int_0^t\sigma(s,y)\dd s}\dd y \]
and $\sigma$ satisfies the second equation of (\ref{eta}).

From the uniqueness of solutions in the Lagrangian formulation, we conclude that if $(\varrho_i, u_i)$, $i=1,2$ are solutions to (\ref{NS2}), then they are equal at the level of Lagrangian coordinates. In particular,
 \[ \ddiv u_1(t,x_1(t,y))=\ddiv u_2(t, x_2(t,y))=\sigma(t,y). \]
Therefore the uniqueness of the solutions to (\ref{NS2}) is equivalent to uniqueness of solutions to the equation
 \begin{equation}\label{div} \ddiv u(t,x(t,y)) = \sigma(t,y), \end{equation}
 where $\sigma\in L^\infty([0,T]\times\mathbb{T}^d)$ is given and $x(t,y)$ satisfies an ODE
 \begin{equation}\label{ode} \begin{aligned} \dot{x}(t,y) &= u(t,x), \\
 x(0,y) &= y. \end{aligned} \end{equation}
 Having the unique $u$, the uniqueness of $\varrho$ follows then again from the classical results from  transport theory. The regularity of $u$ provides that in particular $u\in L^1(0,T;W^{1,p})$ for some $p\geq 1$ and $\ddiv u\in L^1(0,T;L^\infty)$. Therefore there exists a unique $\varrho\in L^\infty([0,T]\times\mathbb{T}^d)$, which is a solution to the continuity equation (\ref{continuity}). 
 
\subsection{Definition of the flow $x_s$}
 For $u_1$ and $u_2$ being the solutions to (\ref{div})-(\ref{ode}), we introduce the \textit{weighted flow between $u_1$ and $u_2$}, that is a function $x_s(t,y), \; s\in[0,1]$ such that $x_s$ for $s=0$ is the flow generated by $u_2$ and for $s=1$ the flow generated by $u_1$. Such $x_s$ is defined by the ordinary differential equation
 
\begin{equation}\label{x_s}
\begin{aligned} \dot{x_s}(t,y) &= s u_1(t,x_s)+(1-s) u_2(t,x_s) \\
x_s(0,y) &= y \end{aligned}\end{equation}
for $s\in [0,1]$. Note that the Jacobian $J_s$ of $x_s$ satisfies the same bounds as the Jacobian $J$ of $x_1$ and $x_2$, namely
\begin{equation}\label{J_s} e^{-L} \leq J_s(t,y) \leq e^{L} \end{equation}
with $\displaystyle L=\int_0^T \|s\ddiv u_1(t,x_s)+(1-s)\ddiv u_2(t,x_s) \|_\infty \dd t \leq \int_0^T \|\sigma\|_\infty \dd t$. \\

The first step to show uniqueness is to obtain certain $L^p$ estimates for the derivative of $x_s$ with respect to $s$.

\begin{lem}\label{reg_xs}
If $x_s$ is defined by (\ref{x_s}) for $u_1,u_2\in L^\infty([0,T]\times\mathbb{T}^d)$ such that $\nabla u_i\in L^\infty(0,T;BMO)$, $i=1,2$, then for sufficiently small $t$, $\frac{\dd x_s}{\dd s}\left(t,\cdot\right)\in L^p$ for some $p>4$.
\end{lem}

Proof: differentiating $\dot{x}_s$ with respect to $s$, we get
\begin{equation}\label{ds} \frac{\dd\dot{x}_s}{\dd s} = u_1(t,x_s)-u_2(t,x_s) + (s\nabla u_1(t,x_s)+(1-s)\nabla u_2(t,x_s))\frac{\dd x_s}{\dd s}, \end{equation}
hence from the Gronwall's lemma
\[ \left|\frac{\dd x_s}{\dd s}\right| \leq \exp\left(\int_0^t |s\nabla u_1(\tau,x_s)+(1-s)\nabla u_2(\tau,x_s)|\dd\tau\right)\int_0^t |u_1(\tau,x_s)-u_2(\tau,x_s)|\dd\tau. \]
Let $\nabla v=s\nabla u_1+(1-s)\nabla u_2$. From the H\"older inequality,
\[ \begin{aligned} \int\left|\frac{\dd x_s}{\dd s}\right|^p\dd y \leq & \left(\int \exp\left(2p\int_0^t |\nabla v(\tau,x_s)|\dd\tau\right)\dd y\right)^{1/2} \\
&\times \left(\int_0^t\int |u_1(\tau,x_s)-u_2(\tau,x_s)|^{2p}\dd y\dd\tau\right)^{1/2}. \end{aligned} \]
As $u_1,u_2\in L^\infty([0,T]\times\mathbb{T}^d)$, the term 
$\displaystyle \left(\int_0^t\int |u_1(\tau,x_s)-u_2(\tau,x_s)|^{2p}\dd y\dd\tau\right)^{1/2} $
is bounded. We will now estimate
\[ I = \int \exp\left(2p\int_0^t |\nabla v(\tau,x_s(\tau,y))\dd\tau\right)\dd y. \]

By Jensen's inequality and the bounds on Jacobian $J_s$, we have
\[ \begin{aligned} \int e^{2p\int_0^t |\nabla v(\tau,x_s(\tau,y))|\dd\tau}\dd y &\leq \int \frac{1}{t}\int_0^t e^{2pt|\nabla v(\tau,x_s(\tau,y))|} \; \dd\tau \dd y \\
&\leq \left\|\frac{1}{J_s}\right\|_{\infty,T}\int \frac{1}{t}\int_0^t e^{2pt|\nabla v(\tau,x_s(\tau,y))|}J_s(\tau,y) \; \dd\tau\dd y \\
&\leq e^L \frac{1}{t}\int_0^t\int e^{2pt|\nabla v(\tau,x)|} \; \dd x \dd\tau. \end{aligned} \]
From the fact that $\sup_{t\in [0,T]}\|\nabla u_i\|_{BMO} \leq C$ for $i=0,1$, we know that $\nabla v\in L^\infty(0,T;BMO)$ and from the Corollary \ref{wniosek_JN}
\[ \int e^{2pt|\nabla v(\tau,x)|} \dd x \leq C \; \; \text{ for all } \; \; p \leq \frac{C}{2t\|\nabla v\|_{L^\infty(0,T;BMO)}}. \]
In particular, for sufficiently small $t$ we have $\frac{\dd x_s}{\dd s}\in L^p(\mathbb{T}^d)$ for some $p>4$. \qed

\begin{lem}\label{reg_xs2}
Let $x_s$ be as in Lemma \ref{reg_xs} and let $T_1$ be such that $\frac{\dd x_s}{\dd s}\in L^p$ for $p>4$ and $t\in[0,T_1]$. Then $\left\|\frac{\dd x_s}{\dd s}\left(t,\cdot\right)\right\|_2^2$ satisfies the inequality
\begin{equation}\label{alpha_s} \frac{\dd}{\dd t}\left\|\frac{\dd x_s}{\dd s}\right\|_2^2 \leq C\left\|\frac{\dd x_s}{\dd s}\right\|_2^2\left(1+\left|\ln\left\|\frac{\dd x_s}{\dd s}\right\|_2^2\right|\right) + C\|u_1(t,\cdot)-u_2(t,\cdot)\|_2^2 \end{equation}
for $t\in [0,T_1]$.
\end{lem}
Proof: multiplying both sides of (\ref{ds}) by $\frac{\dd x_s}{\dd s}$, we obtain
\begin{equation}\label{ds^2} \frac{1}{2}\frac{\dd}{\dd t}\left|\frac{\dd x_s}{\dd s}\right|^2 = \frac{\dd x_s}{\dd s}\nabla v(t,x_s) \frac{\dd x_s}{\dd s} + \big(u_1(t,x_s)-u_2(t,x_s)\big)\frac{\dd x_s}{\dd s}. \end{equation}

Integrating (\ref{ds^2}) over torus, we get
\[\begin{aligned} \frac{\dd}{\dd t}\int \left|\frac{\dd x_s}{\dd s}\right|^2\dd y &\leq 2\left|\int \frac{\dd x_s}{\dd s}\nabla v(t,x_s) \frac{\dd x_s}{\dd s} \dd y \right| + 2\int |u_1(t,x_s)-u_2(t,x_s)|\left|\frac{\dd x_s}{\dd s}\right| \dd y \\ 
&\leq 2\left|\int \frac{\dd x_s}{\dd s}\nabla v(t,x_s) \frac{\dd x_s}{\dd s} \dd y \right| + \int |u_1(t,x_s)-u_2(t,x_s)|^2\dd y + \int \left|\frac{\dd x_s}{\dd s}\right|^2\dd y. \end{aligned}\]
From the bound (\ref{J_s}) on $J_s$, we have 
\[ \int |u_1(t,x_s)-u_2(t,x_s)|^2\dd y \leq C\|u_1(t,\cdot)-u_2(t,\cdot)\|_2^2. \]

By the regularity of $\frac{\dd x_s}{\dd s}$ from Lemma \ref{reg_xs}, $\left|\frac{\dd x_s}{\dd s}\right|^2\in L^q$ for some $q>2$. Therefore we can apply Corollary \ref{wniosek} to estimate $\displaystyle \left|\int \frac{\dd x_s}{\dd s}\nabla v(t,x_s) \frac{\dd x_s}{\dd s} \dd y \right|$. In consequence we obtain the inequality (\ref{alpha_s}), where $C$ depends on $\|\nabla v\|_{BMO}$, $\left\|\frac{\dd x_s}{\dd s}\right\|_q$ and $\|J\|_\infty$. \qed

\subsection{The final argument}
Having the results from the previous subsection, we can now prove the uniqueness:
\begin{tw}\label{uniqueness}
The solution to system (\ref{div})-(\ref{ode}) with regularity from Theorem \ref{main} is unique.
\end{tw} 
Proof: Let $u_1,u_2$ be the solutions to (\ref{div})-(\ref{ode}), $u_i=\nabla\phi_i$ and $x_1,x_2$ are the flows generated by $u_1,u_2$ respectively. We will show that $\|u_1-u_2\|_2=0$ for all $t\in [0,T]$. By the weak formulation of $\ddiv u(t,x(t,y))=\sigma(t,y)$, for any $\xi\in C^
\infty([0,T]\times\mathbb{T}^d)$ we have
\[\begin{aligned} \int (u_1(t,x)-u_2(t,x))\nabla\xi(t,x) \dd x &= -\int (\ddiv u_1(t,x)-\ddiv u_2(t,x))\xi(t,x) \dd x \\
&= -\int \sigma(t,y)J(t,y)(\xi(t,x_1(t,y))-\xi(t,x_2(t,y)))\dd y. \end{aligned} \]
Using definition of the flow $x_s$, we can rewrite the last integral as
\begin{multline*} \int \sigma(t,y)J(t,y)\int_0^1 \frac{\dd}{\dd s}\xi(t,x_s(t,y))\dd s\dd y 
= \int \sigma(t,y)J(t,y)\int_0^1 \nabla\xi(t,x_s(t,y))\frac{\dd x_s(t,y)}{\dd s}\dd s\dd y. \end{multline*}
By the density of smooth functions in $W^{1,2}$, we can choose $\xi=\phi_1-\phi_2$. Then $\nabla\xi=u_1-u_2$ and we obtain
\[\begin{aligned} \int |u_1-u_2|^2\dd x &= -\int \sigma(t,y)J(t,y)\int_0^1 (u_1(t,x_s(t,y))-u_2(t,x_s(t,y)))\frac{\dd x_s}{\dd s}(t,y) \dd s\dd y \\
&\leq \|\sigma\|_\infty\|J\|_\infty \|u_1-u_2\|_2 \int_0^1 \left\|\frac{\dd x_s}{\dd s}\right\|_2 \dd s. \end{aligned}\]
Hence
\begin{equation}\label{u1-u2} \|u_1-u_2\|_2 \leq C\int_0^1 \left\|\frac{\dd x_s}{\dd s}\right\|_2 \dd s. \end{equation}
 Substituting (\ref{u1-u2}) into (\ref{alpha_s}) and integrating over $s$, we get
\[ \frac{\dd}{\dd t}\int_0^1 \left\|\frac{\dd x_s}{\dd s}\right\|_2^2 \; \dd s \leq C\int_0^1 \left\|\frac{\dd x_s}{\dd s}\right\|_2^2\left(1+\left|\ln\left\|\frac{\dd x_s}{\dd s}\right\|_2^2\right|\right) \; \dd s + C\left(\int_0^1 \left\|\frac{\dd x_s}{\dd s}\right\|_2 \dd s\right)^2 \]
for $t\in [0,T_1]$ for some $T_1\leq T$. \\
Let $\alpha(t)=\int_0^1 \left\|\frac{\dd x_s}{\dd s}\right\|_2^2 \dd s$. As the function $x(1+|\ln x|)$ is concave for $x<1$ and $x^2$ is convex, we can estimate both terms in the right hand side from Jensen's inequality and obtain
\[ \dot{\alpha} \leq C\alpha(1+|\ln\alpha|)+ C\alpha. \]
From Osgood's lemma, the problem
\[ \dot{z} = Cz(1+|\ln z|), \quad z(0) = 0 \]
has a unique solution $z\equiv 0$. Therefore, as $\frac{\dd x_s}{\dd s}_{|_{t=0}}=0$, we have 
\[ \alpha(t) \leq 0 \textrm{ for all } t\in [0,T_1]. \]
In conclusion, $\frac{\dd x_s}{\dd s}\equiv 0$ for all $t\in [0,T_1]$ and so is $\|u_1-u_2\|_2$. Having that, we can perform analogous reasoning on the consecutive intervals $[nT_1,(n+1)T_1]$ to get $u_1=u_2$ for all $t\in[0,T]$. \qed

 \section{The existence of solutions to (\ref{div})-(\ref{ode})}\label{euler_sect}
 In this section we prove that the transformation from equation (\ref{div}) together with (\ref{ode}) to equation (\ref{NS2}) is well defined, which will end the proof of Theorem \ref{main}. Having the solution $(\eta,\sigma)$ in the Lagrangian coordinates by Theorem \ref{lagrange}, we define the transformation to Eulerian coordinates $(\varrho,u)$. In other words, we need to find $u$ such that $\displaystyle \int u(t,x)\dd x=0$ and
$u$ satisfies (\ref{div})-(\ref{ode}). By virtue of the discussion at the beginning of section \ref{uniq_sect}, such $u$ provides us also existence of the density $\varrho$.
 
 \begin{tw}\label{to_euler} Let $\sigma\in L^\infty([0,T]\times\mathbb{T}^d)$. There exists $u\in L^\infty([0,T]\times\mathbb{T}^d)$ such that $u$ is a solution to system (\ref{div})-(\ref{ode}) and
 \[ \ddiv u\in L^\infty([0,T]\times\mathbb{T}^d), \; \nabla u \in L^\infty(0,T;BMO). \]
 \end{tw}
 
 First, we prove the existence for smoothened $\sigma$, by putting 
 \[ \sigma_\delta=\sigma\ast\kappa_\delta, \]
 where $\kappa_\delta$ is a standard mollifier.
 
 \begin{lem}\label{lemma_sigma}
 There exists a unique $u_\delta\in C(0,T;W^{1,\infty})$ satisfying
 \[ \ddiv u_\delta(t,x_\delta(t,y)) = \sigma_\delta(t,y), \]
 where $x_\delta(t,y)$ is given by (\ref{ode}).
 \end{lem}
 
  Proof: We define the suitable map, and then apply the Banach fixed point theorem. Let 
  \[\Phi\colon C(0,T;W^{1,\infty}) \to C(0,T;W^{1,\infty}) \]
 be defined in the following way:
 \begin{enumerate}
     \item If $\bar{u}\in C(0,T;W^{1,\infty})$, then $\bar{u}$ is Lipschitz, so there exists a unique solution to system
     \begin{equation}\label{rrz2}
       \dot{x} = \bar{u}(t,x), \quad x(0) = y.  \end{equation}
     \item We can now invert $x(t,y)$ to get $y(t,x)$ instead. After differentiation of (\ref{rrz2}) with respect to $y$, we get an ODE for the matrix $H(t,y)=\frac{\partial x}{\partial y}(t,y)$:
     \begin{equation} \partial_t H = \nabla_x\bar{u}(t,x(t,y))H, \; H(0,y)=I. \end{equation}
     Moreover, the equation for $J(t,y)=\det H(t,y)$ yields
     \begin{equation}\label{J} \partial_tJ(t,y) = \text{div}_x\bar{u}(t,x(t,y))J(t,y). \end{equation}
     Therefore we have the estimates
     \[ \exp\left( -\int_0^T\| \nabla\bar{u}\|_\infty \dd s\right) \leq \| H\|_\infty \leq \exp\left(\int_0^T\|\nabla\bar{u}\|_\infty\dd s\right), \]
     and
     \[ \exp\left( -\int_0^T\| \text{div}\bar{u}\|_\infty \dd s\right) \leq \| J\|_\infty \leq \exp\left(\int_0^T\|\text{div}\bar{u}\|_\infty\dd s\right), \]
     and $H$ is invertible, which allows us to treat $y$ as a function of $x$.
     \item Finally, we put $\Phi(\bar{u})=u$, where $u$ is a unique solution to the system
     \begin{equation}\label{eliptyczne}
     \begin{aligned} u(t,x) &=\nabla\phi(t,x), \\
     \Delta\phi(t,x) &= \sigma_\delta(t,y(t,x)) \end{aligned} \end{equation}
for $y(t,x)$ being the inverse flow associated with $\bar{u}$. \\

     The function $\sigma_\delta(t,y(t,\cdot))$ is Lipschitz, as
     \[\begin{aligned} |\sigma_\delta(t,y(t,x_1))-\sigma_\delta(t,y(t,x_2))| &\leq \|\nabla\sigma_\delta\|_\infty |y(t,x_1)-y(t,x_2)| \\
     &\leq \|\nabla\sigma_\delta\|_\infty\left\|\frac{\dd y}{\dd x}\right\|_\infty |x_1-x_2|. \end{aligned}\]
 Therefore $\phi(t,\cdot)\in C^{2,1}(\mathbb{T}^d)$ and in consequence $u(t,\cdot)\in C^{1,1}(\mathbb{T}^d)$ and we have the estimates
 \begin{equation}\label{holder}
 \begin{aligned}
     \sup_{0\leq t\leq T}\|u\|_\infty &\leq C\|\sigma_\delta\|_{\infty,T}, \\
     \sup_{0\leq t\leq T}\|\nabla u\|_{C^{0,1}} &\leq C\sup_{0\leq t\leq T}\|\sigma_\delta(t,y(t,\cdot))\|_{C^{0,1}}. \end{aligned}
 \end{equation}
 
 Moreover, as $\sigma\in L^\infty([0,T]\times\mathbb{T}^d)$, by the classical results for $L^p$ regularity of strong solutions to Poisson equation (see e.g. \cite{wu_yin_wang}), for any $p\in (1,\infty)$ we have the estimate
 \begin{equation}\label{oszac}
 \sup_{0\leq t\leq T}\|u\|_\infty + \sup_{0\leq t\leq T} \|u\|_{W^{1,p}} + \sup_{0\leq t\leq T}\|\nabla u\|_{BMO} \leq C\|\sigma\|_{\infty,T},
 \end{equation}
 which are uniform with respect to $\delta$.
 \end{enumerate}
 
 \subsection{Additional regularity of $u$} 
We will now show that the fixed points of $\Phi$ are uniformly bounded in $L^\infty(0,T;W^{2,p})$. By use of an appropriate logarithmic inequality, it implies that $\Phi(K)\subseteq K$ for some bounded and closed $K\subset C(0,T;W^{1,\infty})$.
 \begin{prop}
 If $u_\delta$ is a fixed point of $\Phi$, then
 \begin{equation}\label{1,p} \sup_{0\leq t\leq T}\|\nabla u_\delta\|_{W^{1,p}} \leq C(p,\delta,T). \end{equation}
 \end{prop}
 Proof: After differentiation of (\ref{eliptyczne}) with respect to $x_i$, we get
 \[ \Delta_x\frac{\partial\phi}{\partial x_i} = \nabla_y\sigma_\delta\cdot\frac{\partial y}{\partial x_i}. \]
 Hence for any $1<p<\infty$ the standard elliptic estimate gives
 \[ \|\nabla\phi\|_{W^{2,p}} \leq C(p)\|\nabla\sigma_\delta\|_\infty\left\|\frac{\dd y}{\dd x}\right\|_\infty, \]
 which leads to
 \begin{equation}\label{2,p} \sup_{0\leq t\leq T}\|\nabla u\|_{W^{1,p}} \leq C(p)\|\nabla\sigma_\delta\|_{\infty,T}e^{T\|\nabla\bar{u}\|_{\infty,T}}. \end{equation}
 Now take $\bar{u}=u=u_\delta$ in (\ref{2,p}) and apply the following inequality for $p>d$:
 \begin{equation}\label{bmo_est}
     \|\nabla f\|_\infty \leq C(p)\left(1+\|\nabla f\|_{BMO}\left(\ln^+(\|\nabla f\|_{W^{1,p}}+\|f\|_\infty)\right)^{1/2}\right).
 \end{equation}
 The proof of the above estimate one can find in \cite{bmo_ineq}, Corollary 2.4. By the estimates (\ref{oszac}) and Cauchy inequality, we have
 \[\begin{aligned} \|\nabla u_\delta\|_{L^\infty(0,T;W^{1,p})} &\leq C(p)\|\nabla\sigma_\delta\|_{\infty,T}\exp\left(T\sup_{0\leq t\leq T}\|\nabla u_\delta\|_{\infty}\right) \\
 &\leq \exp\left(CT\sup_{0\leq t\leq T}\|\nabla u_\delta\|_{BMO}\left(\ln^+(\|\nabla u_\delta\|_{L^\infty(0,T;W^{1,p})}+\|u_\delta\|_{\infty,T})\right)^{1/2}\right) \\
 &\leq \exp\left(CT^2+\frac{1}{2}\ln^+\left(\|\nabla u_\delta\|_{L^\infty(0,T;W^{1,p})}+\|u_\delta\|_{\infty,T}\right)\right) \\
 &\leq C\left(\|\nabla u_\delta\|_{L^\infty(0,T;W^{1,p})}+\|u_\delta\|_{\infty,T}\right)^{1/2}.
 \end{aligned}\]
 Therefore
 \[ \|\nabla u_\delta\|_{L^\infty(0,T;W^{1,p})}^2 - C\|\nabla u_\delta\|_{L^\infty(0,T;W^{1,p})} - C \leq 0 \]
 and in consequence 
 \[ \|\nabla u_\delta\|_{L^\infty(0,T;W^{1,p})} \leq C \]
 for some $C$ depending on $p$, $\delta$ and $T$. \qed \\
 
The analogous reasoning also provides that $\Phi(K_{T_1})\subseteq K_{T_1}$ for
\[ K_{T_1} = \{u\in C(0,T_1;W^{1,\infty}): \ \|u\|_{\infty,T_1}, \|\nabla u\|_{L^\infty(0,T_1;BMO)}\leq C_1 \ \text{ and } \ \|\nabla u\|_{\infty,T_1} \leq C_2 \}, \]
where $C_1$ is the constant from estimates (\ref{oszac}) and $C_2=C(p)\left(1+C_1\left(\ln^+(C_p+C_1)\right)^{1/2}\right)$ is the right hand side of (\ref{bmo_est}) for $C_p$ being the constant from (\ref{1,p}).
 
\subsection{The fixed point argument.}
\paragraph{Local existence.} We show that $\Phi\colon K\to K$ is a contraction for sufficiently small $T_1$. From the elliptic estimates, we have
\[\begin{aligned} \|u_1(t,\cdot)-u_2(t,\cdot)\|_\infty &\leq C\|\sigma_\delta(t,y_1(t,\cdot))-\sigma_\delta(t,y_2(t,\cdot))\|_\infty \\
&\leq C\|\nabla\sigma_\delta\|_{\infty,T}\|y_1(t,\cdot)-y_2(t,\cdot)\|_\infty. \end{aligned}\]
By Lemma \ref{lemat} from the Appendix, for small $t$ we have
\[\|y_1(t,\cdot)-y_2(t,\cdot)\|_\infty \leq Ct\|\bar{u}_1-\bar{u}_2\|_{\infty,T} \]
 and due to the uniform bound on $\|\nabla u\|_\infty$ for $u\in K$ the constant $C$ does not depend on $\bar{u}_1$ and $\bar{u}_2$. Hence
\begin{equation}\label{z_lematu} \|u_1-u_2\|_{\infty,T_1} \leq C\|\nabla\sigma_\delta\|_{\infty,T}(T_1+o(T_1))\|\bar{u}_1-\bar{u}_2\|_{\infty,T_1}. \end{equation}
 Choosing $T_1$ such that $C\|\nabla\sigma_\delta\|_{\infty,T}T_1 < 1$, we get that $\Phi$ is a contraction on $K$. Therefore there exists a unique fixed point $u_\delta$ of $\Phi$ on the interval $[0,T_1]$. 
 
 \paragraph{Extension to $[0,T]$.} Having the uniqueness on $[0,T_1]$, we are able to perform the same reasoning on $[T_1,2T_1]$. Note that the estimate (\ref{z_lematu}) on $[T_1,2T_1]$ would still depend only on the length of the interval. Therefore again $\Phi$ is a contraction on $[T_1,2T_1]$, which gives us the unique fixed point on $[0,2T_1]$. Performing this procedure on the consecutive intervals $[nT_1,(n+1)T_1]$, we obtain the existence of a unique fixed point on the whole interval $[0,T]$, which completes the proof of Lemma \ref{lemma_sigma}. \qed

 \subsection{Letting $\delta\to 0$.}
 In Lemma \ref{lemma_sigma} we obtained the unique $u_\delta=\nabla\phi_\delta$, which satisfies the equation
 \begin{equation}\label{phi_eq} \Delta\phi_\delta(t,x) = \sigma_\delta(t, y_\delta(t,x)) \end{equation}
 and the uniform estimates (\ref{oszac}). Moreover, it turns out that $\partial_t u_\delta$ is uniformly bounded in $L^2$. Indeed, from the weak formulation of (\ref{phi_eq}), for any $\pi\in C^\infty(\mathbb{T}^d)$ we have
 \[ \int \Delta\phi_\delta(t,x)\pi(x) \; \dd x = \int \sigma_\delta(t,y_\delta(t,x))\pi(x) \; \dd x = \int \sigma_\delta(t,y)\pi(x_\delta(t,y))J_\delta(t,y)\; \dd y,\]
 where $J_\delta(t,y)=\exp\left(\int_0^t\sigma_\delta(s,y)\dd s\right)$ is the Jacobian of $x_\delta$. Differentiating this equality with respect to time, we obtain
 \begin{equation}\label{phi_t} \begin{aligned}
 \int \Delta\partial_t\phi_\delta(t,x)\pi(x) \; \dd x = & \int \partial_t\sigma_\delta(t,y)\pi(x_\delta(t,y))J_\delta(t,y) \; \dd y \\
 &+ \int \sigma_\delta(t,y)\partial_t\pi(x_\delta(t,y))J_\delta(t,y) \; \dd y \\
 &+ \int \sigma_\delta(t,y)\pi(x_\delta(t,y))\partial_tJ_\delta(t,y) \; \dd y.
 \end{aligned} \end{equation}
 Let us now estimate the terms on the right hand side of (\ref{phi_t}). First, observe that from equation (\ref{eta}) $\partial_t\sigma$ is bounded, and therefore 
 \[ |\partial_t\sigma_\delta|=|(\partial_t\sigma)\ast\kappa_\delta|\leq |\partial_t\sigma| \in L^\infty([0,T]\times\mathbb{T}^d). \]
 For the second term, we have 
 \[ \partial_t\pi(x_\delta(t,y))=\nabla\pi(x_\delta(t,y))\dot{x_\delta}(t,y) = \nabla\pi(x_\delta(t,y))u(t,x_\delta(t,y)) \]
 and then using uniform estimates on $u_\delta$, we get
 \[ \int |\nabla\pi(x_\delta(t,y))|^2|u_\delta(t,x_\delta(t,y))|^2 \; \dd y \leq \left\|\frac{1}{J_\delta}\right\|_\infty\int |\nabla\pi(t,x)|^2|u_\delta(t,x)|^2\; \dd x \leq C, \]
 therefore $\partial_t\pi(x_\delta(t,y))\in L^\infty(0,T;L^2)$. The third term is bounded as well, as $\partial_t J_\delta = \sigma_\delta J_\delta$
 and both $\sigma_\delta$ and $J_\delta$ are bounded by some $C(\|\sigma\|_{\infty,T})$. \\
 The above estimates imply that $\Delta\partial_t\phi_\delta$ is bounded in $L^\infty(0,T;W^{-1,2})$ uniformly in $\delta$ and therefore 
 \begin{equation}\label{u_t}
     \|\partial_t u_\delta\|_{L^\infty(0,T;L^2)} = \|\nabla\partial_t\phi_\delta\|_{L^\infty(0,T;L^2)} \leq C.
 \end{equation} \\
 
 We now let $\delta\to 0^+$ and therefore obtain the solution to equation (\ref{div}). The estimates (\ref{oszac}) give
 \[ \|\phi_\delta\|_{L^\infty(0,T;W^{2,p})} \leq C, \]
 so $\phi_\delta\rightharpoonup ^* \phi$ in $L^\infty(0,T;W^{2,p})$ up to a subsequence. From the uniform estimates (\ref{oszac}) and (\ref{u_t}), Aubin-Lions Lemma implies that in particular $u_\delta$ is compact in $L^1([0,T]\times\mathbb{T}^d)$. Moreover, using Theorem 2.9 from \cite{crippa-delellis}, we get
 \[ \sup_{0\leq t\leq T}\|x(t,y)-x_\delta(t,y)\|_{L^1(\mathbb{T}^d)} \leq C\left|\ln\left(\|u-u_\delta\|_{L^1([0,T]\times\mathbb{T}^d)}\right)\right|^{-1}, \]
 where $x(t,y)$ is the flow generated by this weak* limit $u=\nabla\phi$. Therefore $x_\delta\to x$ in $L^\infty(0,T;L^1)$. \\
 
 The above convergence allows us to pass to the limit with $\delta\to 0$ in a weak formulation of (\ref{eliptyczne}). For any $\xi\in C^\infty([0,T]\times\mathbb{T}^d)$, we have
 \[\begin{aligned} \int_0^T\int \Delta\phi_\delta(t,x)\xi(t,x)\dd x \dd t &= \int_0^T\int (\sigma\ast\kappa_\delta)(t,y_\delta(t,x))\xi(t,x)\dd x \dd t \\
 &= \int_0^T\int (\sigma\ast\kappa_\delta)(t,y)\xi(t,x_\delta(t,y))J_\delta(t,y)\dd y \dd t. \end{aligned}\]
 Letting $\delta\to 0$, we get
 \begin{equation}\label{limit} \int_0^T\int \Delta\phi(t,x)\xi(t,x)\dd x \dd t = \int_0^T\int\sigma(t,y)\xi(t,x(t,y))J(t,y)\dd y \dd t, \end{equation}
 where $J=\exp\left(\int_0^t \sigma(s,y)\dd s\right)$ is the Jacobian of the limit flow $x(t,y)$.
 
 To deduce that indeed we have $\ddiv u(t,x(t,y))=\sigma(t,y)$, we need to change the variables in one of the sides in (\ref{limit}). Despite the fact that $x(t,\cdot)$ is not a diffeomorphism, Lemma 3.1 from \cite{colombo-crippa} allows us to perform the change of variables in the left hand side of (\ref{limit}) and obtain
 \[ \int_0^T\int \ddiv u(t,x(t,y)) \xi(t,x(t,y)) J(t,y)\dd y \dd t. \]
 Therefore the equality (\ref{limit}) is transformed into
 \[ \int_0^T\int \big[\ddiv u(t,x(t,y))-\sigma(t,y)\big]\xi(t,x(t,y))J(t,y)\dd y \dd t = 0. \]
 As $\ddiv u\in L^\infty([0,T]\times\mathbb{T}^d)$, the Jacobian $J(t,y)$ is strictly positive. Hence, from the arbitrary choice of $\xi$, we have $\ddiv u(t,x(t,y))=\sigma(t,y)$ in the sense of distributions, which ends the proof of Theorem \ref{to_euler}. \qed
 
As by Theorem \ref{lagrange} the norms $\|\eta\|_{\infty,T}$ and $\|\sigma\|_{\infty,T}$ does not depend on $T$, so are $\|\varrho\|_{\infty,T}$, $\|\ddiv u\|_{\infty,T}$ and the estimates given by (\ref{oszac}). Then again from arbitrary choice of $T$ we obtain the unique existence on the whole real line and hence the proof of Theorem \ref{main} is completed.
 
 \paragraph{Acknowledgments.}  The author has been supported by National Science Centre grant \\ 2018/29/B/ST1/00339 (Opus) and would like to thank Piotr Mucha and Tomasz Piasecki for their invaluable help and remarks through the process of creating the paper.
 
 \begin{appendices}
 
 \section{Existence of solutions to (\ref{eta})}\label{existence_proof}
 We here prove that there exists a unique solution to (\ref{eta}), which completes the proof of Theorem \ref{lagrange}.
 \begin{lem}\label{ex_L}
  For $\varrho_0\in L^\infty(\mathbb{T}^d)$ and any $T>0$ there exists a unique global in time solution 
  \[(\eta, \sigma)\in L^\infty([0,T]\times\mathbb{T}^d)\times L^\infty([0,T]\times\mathbb{T}^d) \]
  to the equation (\ref{eta}) with the initial condition $\eta(0,y)=\varrho_0(y)$.
 \end{lem}
 Proof: The proof relies on double application of the Banach fixed point theorem. First, observe that for a fixed $\eta\in L^\infty([0,T]\times\mathbb{T}^d)$, there exists a unique $\sigma$, satisfying
 \begin{equation}\label{sigma} \sigma = p(\eta)-\frac{1}{|\mathbb{T}^d|}\int p(\eta(t,y))\exp\left( \int_0^t \sigma(s,y)\dd s\right) \dd y. \end{equation}
 To see that, take $\Psi\colon L^\infty([0,T]\times\mathbb{T}^d)\to L^\infty([0,T]\times\mathbb{T}^d)$ such that
 \[ \Psi(\sigma)=p(\eta)-\frac{1}{|\mathbb{T}^d|}\int p(\eta(t,y))\exp\left( \int_0^t \sigma(s,y)\dd s\right) \dd y. \]
 Then 
 \[\begin{aligned} |\Psi(\sigma_1)-\Psi(\sigma_2)| &\leq \frac{1}{|\mathbb{T}^d|}\int |p(\eta(t,y))|\left|\exp\left( \int_0^t \sigma_1(s,y)\dd s\right)-\exp\left(\int_0^t \sigma_2(s,y)\dd s\right)\right| \dd y. \\
 &\leq \frac{1}{|\mathbb{T}^d|}\int |p(\eta(t,y)|\int_0^t |\sigma_1(s,y)-\sigma_2(s,y)|\dd s \dd y \\
 &\leq Ct\|\sigma_1-\sigma_2\|_{\infty,\tau}.
 \end{aligned}\]
 Hence taking $\tau$ such that $C\tau<1$, from the Banach fixed point theorem we get the existence of a unique solution on the interval $[0,\tau]$. However, as $\sigma$ does not blow up on an interval $[0,T]$, we are able to extend the solution to (\ref{sigma}) to the whole interval. \\
 
 Now define $\Phi\colon L^\infty([0,T]\times\mathbb{T}^d)\to L^\infty([0,T]\times\mathbb{T}^d)$ as
 \[ \Phi(\eta) =\varrho_0(y)-\int_0^t\eta(s,y)\sigma(s,y)\mathrm{d}s, \]
 where $\sigma$ is given by (\ref{sigma}). Similarly as before, we want to obtain the estimate
 \[ \|\Phi(\eta_1)-\Phi(\eta_2)\|_{\infty,\tau} \leq C\tau\|\eta_1-\eta_2\|_{\infty,\tau} \]
 and choose $\tau$ such that $C\tau<1$ and $\Phi\colon \{\eta: \|\eta\|_{\infty,\tau}< r\}\to\{\eta: \|\eta\|_{\infty,\tau}< r\}$ is a contraction. We have
 \[ |\Phi(\eta_1)-\Phi(\eta_2)|=|\int_0^t\eta_1\sigma_1-\eta_2\sigma_2 \mathrm{d}s| \leq \int_0^t|\eta_1-\eta_2\|\sigma_1|\mathrm{d}s+\int_0^t\eta_2|\sigma_1-\sigma_2|\mathrm{d}s. \]
The first integral can be bounded by $Ct\|\eta_1-\eta_2\|_{\infty,\tau}$, so to complete the desired estimate we need to estimate the difference of $\sigma_1$ and $\sigma_2$. We have
 \[ |\sigma_1-\sigma_2|\leq |p(\eta_1)-p(\eta_2)|+|\{p(\eta_1)\}_{\sigma_1}-\{p(\eta_2)\}_{\sigma_2}|. \]
 As the derivative of $p$ is bounded on an interval $[0,r]$, we can estimate the first element by $C\|\eta_1-\eta_2\|_\infty$. Moreover,
 \[\begin{aligned} |\mathbb{T^d}| & |\{p(\eta_1)\}_{\sigma_1}-\{p(\eta_2)\}_{\sigma_2}| = \left|\int p(\eta_1)\exp\left(\int_0^t\sigma_1\mathrm{d}s\right)-p(\eta_2)\exp\left(\int_0^t\sigma_2\mathrm{d}s\right)\mathrm{d}y \right| \\
&\leq \int \exp\left(\int_0^t\sigma_1\mathrm{d}s\right)|p(\eta_1)-p(\eta_2)|\mathrm{d}y + \int p(\eta_2)\left|\exp\left(\int_0^t\sigma_1\mathrm{d}s\right)-\exp\left(\int_0^t\sigma_2\mathrm{d}s\right)\right|\mathrm{d}y \\
&\leq C\|\eta_1-\eta_2\|_\infty + C\int\left| \int_0^t(\sigma_1-\sigma_2) \mathrm{d}s\right| \mathrm{d}y \\
&\leq C\|\eta_1-\eta_2\|_\infty + Ct\sup_{0\leq s\leq t}|\{p(\eta_1)\}_{\sigma_1}-\{p(\eta_2)\}_{\sigma_2}|.       \end{aligned}\]
Hence for $\tau$ small enough, we get
\[ \sup_{0\leq t\leq\tau} |\{p(\eta_1)\}_{\sigma_1}-\{p(\eta_2)\}_{\sigma_2}| \leq C\sup_{0\leq t\leq\tau}\|\eta_1-\eta_2\|_\infty, \]
which gives us the desired estimate and ends the proof. \qed
 
 \section{Estimate for the inverse flows.}

 In this section we prove the useful lemma for estimating the difference of the inverse flows $y_i=x_i^{-1}(t,\cdot)$ by the difference of vector fields generating $x_1,x_2$. 
 \begin{lem}\label{lemat} Consider two ordinary differential equations with the same initial value:
 \[ \begin{aligned} \dot{x}_1 &=u_1(t,x_1), \\
\dot{x}_2 &= u_2(t,x_2), \\
x_1(0) &= x_2(0) = y, \end{aligned} \]
where $u_1,u_2\in C(0,T;W^{1,\infty})$. Let $y_1(t,x)$ and $y_2(t,x)$ be the inversions of $x_1$ and $x_2$ with respect to $y$. Then for sufficiently small $t$
\[ \|y_1(t,\cdot)-y_2(t,\cdot)\|_\infty \leq Ct\|u_1-u_2\|_{\infty,T}, \]
where $C=C\left(||\nabla u_1\|_{\infty,T}\right)$.
 \end{lem}
 Proof: Let $M=\|\nabla u_1\|_{\infty,T}$. We have
 \[\begin{aligned} |u_1(t,x_1)-u_2(t,x_2)| &\leq |u_1(t,x_1)-u_1(t,x_2)|+|u_1(t,x_2)-u_2(t,x_2)| \\
&\leq \|\nabla u_1\|_\infty |x_1-x_2| + \| u_1-u_2\|_\infty.
\end{aligned} \]
Substituting it into the difference of $x_1$ and $x_2$, we get
\[\begin{aligned} 
|x_1-x_2| &\leq \int_0^t |u_1(s,x_1)-u_2(s,x_2)|\mathrm{d}s \\
 &\leq \int_0^t \|\nabla u_1\|_\infty |x_1-x_2| \; \mathrm{d}s + \int_0^t \|u_1-u_2\|_\infty \mathrm{d}s.
\end{aligned} \]
Hence from the Gronwall's lemma,
\[\begin{aligned} \|x_1-x_2\|_\infty &\leq \int_0^t \|u_1-u_2\|_\infty \mathrm{d}s + \int_0^t\|\nabla u_1\|_\infty \exp\left(\int_s^t\|\nabla u_1\|_\infty \mathrm{d}\tau\right)\int_0^s \|u_1-u_2\|_\infty \mathrm{d}\tau \mathrm{d}s \\
&\leq \|u_1-u_2\|_{\infty,T}\int_0^t 1+Me^{M(t-s)}s \; \dd s \\
 &\leq Ct\|u_1-u_2\|_{\infty,T} \end{aligned}\]
Analogously, we obtain the estimate
\[\begin{aligned} |x_1(t,y_1)-x_2(t,y_2)| &\leq |x_1(t,y_1)-x_1(t,y_2)|+|x_1(t,y_2)-x_2(t,y_2)| \\
&\leq \left\|\frac{\partial x_1}{\partial y}\right\|_\infty |y_1-y_2|+\|x_1-x_2\|_\infty \\
&\leq \exp\left(\int_0^t\|\nabla u_1\|_\infty\mathrm{d}s\right)|y_1-y_2|+Ct\|u_1-u_2\|_{\infty,T}.
\end{aligned} \]
Combining the above estimates, we get
\begin{align*} |y_1(t,x)-y_2(t,x)| &\leq \int_0^t |u_1(s,x_1(s,y_1(t,x)))-u_2(s,x_2(s,y_2(t,x)))| \mathrm{d}s \\
&\leq \int_0^t \|\nabla u_1\|_\infty |x_1(s,y_1(t,x))-x_2(s,y_2(t,x))| + \|u_1-u_2\|_\infty \mathrm{d}s \\
 &\leq \int_0^t \|\nabla u_1\|_\infty\exp\left(\int_0^s\|\nabla u_1\|_\infty \mathrm{d}\tau\right)|y_1(t,x)-y_2(t,x)|\mathrm{d}s \\
 & \ \ \  + \int_0^t \|\nabla u_1\|_\infty Cs\|u_1-u_2\|_{\infty,T}\mathrm{d}s + \int_0^t \|u_1-u_2\|_\infty \dd s  \\
&\leq \left( e^{Mt}-1\right)|y_1(t,x)-y_2(t,x)| + (t+Ct^2)\|u_1-u_2\|_{L^\infty([0,T]\times\mathbb{T}^d)}.
\end{align*} 
For small $t$ we have $e^{Mt}<2$ and therefore
\[ |y_1(t,x)-y_2(t,x)| \leq \frac{t+Ct^2}{2-e^{Mt}}\|u_1-u_2\|_{\infty,T} \leq Ct\|u_1-u_2\|_{\infty,T} \]

\section{Properties of the $BMO$ space}
We present here the useful properties of the $BMO$ functions, which can be found for example in \cite{stein1} and \cite{torchinsky}.
\begin{de}
 A function $f\in L^1(\Omega)$ belongs to space of bounded mean oscillation $BMO(\Omega)$ iff
 \[ \|f\|_{BMO} = \sup_{Q\subset\Omega}\frac{1}{|Q|}\int_Q |f-\{f\}_Q|\dd x < \infty, \]
 where the supremum is taken over all cubes in $\Omega$.
\end{de}
Note that $\|\cdot\|_{BMO}$ is not a norm, as $\|f\|_{BMO}=0$ for $f$ constant. However, we can equip the space $BMO(\Omega)$ with the norm
\[ \|\cdot\|_{L^1} + \|\cdot\|_{BMO} \]
and then it becomes the Banach space. 

It is straightforward from the definition that the standard mollification is bounded in $BMO$:
\begin{prop}\label{prop_bmo}
For $f\in BMO$ and $\kappa_\delta$ the standard mollifier we have
\[ \|f\ast\kappa_\delta\|_{BMO} \leq \|f\|_{BMO}. \]
\end{prop} 

One of the important tools concerning the $BMO$ spaces is the John-Nirenberg inequality:

\begin{lem}[John-Nirenberg]
There exist constants $c_1,c_2>0$ such that for any cube $Q\subset\Omega$ and $f\in BMO(\Omega)$
\[ |\{x\in Q: |f-\{f\}_Q|>\lambda\}| \leq c_1\exp\left(-\frac{c_2\lambda}{\|f\|_{BMO}}\right)|Q|.  \]
\end{lem}

The useful applications of the John-Nirenberg inequality are the following:
\begin{cor}\label{wniosek_JN}
 Let $f\in BMO(\Omega)$. Then 
 \begin{enumerate}
 \item $f\in L_{loc}^p(\Omega)$ for any $1\leq p<\infty$.
 \item $\displaystyle \sup_{Q\subset\Omega}\int_Q \exp\left(\frac{|f-\{f\}_Q|}{\|f\|_{BMO}}\right) \dd x <\infty. $
 \end{enumerate}
\end{cor}

\subsection{The logarithmic inequality}
 We recall here the inequality from \cite{mucha-rusin}:
 \begin{lem}
Let $f\in BMO(\mathbb{R}^d)$ with compact support and $g\in L^1(\mathbb{R}^d)\cap L^\infty(\mathbb{R}^d)$. Then
\[ \left| \int_{\mathbb{R}^d}fg\dd x\right| \leq C\|f\|_{BMO}\|g\|_{L^1}\big(|\ln\|g\|_{L^1}| + \ln(e+\|g\|_{L^\infty}) \big). \]
\end{lem}
It turns out that after slight modifications, the similar inequality holds for $g\in L^q(\mathbb{T}^d)$ for sufficiently large $q$.
\begin{lem}\label{log_ineq}
Let $f\in BMO(\mathbb{T}^d)$ and $g\in L^q(\mathbb{T}^d)$ for some $q>2$. Then
\begin{equation}\label{mr_q}
\begin{aligned} \left| \int_{\mathbb{T}^d}fg\dd x\right| \leq & C\|f\|_{BMO}\|g\|_{L^1} \\
& \times \big(|\ln\|g\|_{L^1}| + \ln(e+\|g\|_{L^q}) + (1+|\ln\|g\|_{L^1}|)\|g\|_{L^q}^{\frac{q-2}{2}}\big). \end{aligned} \end{equation}
\end{lem}
Proof: Assume $\int g\dd x = 0$. Then $g\in\mathcal{H}^1$ and from duality of $\mathcal{H}^1$ and $BMO$ we have
\[ \left|\int fg \; \dd x \right| \leq \|f\|_{BMO}\|g\|_{\mathcal{H}^1}. \]
By the characterization of $\mathcal{H}^1$ (see e.g. paragraph III.4 in \cite{stein1}) we can write $\|g\|_{\mathcal{H}^1}$ as 
\[ \|g\|_{\mathcal{H}^1} = \|g\|_{L^1} + \sum_{k=1}^d \|R_k g\|_{L^1}, \]
where $R_k$ is the Riesz transform given as $\mathcal{F}(R_k g) = -i\frac{\xi_k}{|\xi|}\mathcal{F}(g)$. As the Riesz transform is an operator of weak-type $(1,1)$ and strong-type $(p,p)$, we can apply Proposition V.3.2. from \cite{torchinsky} and obtain
\begin{equation}\label{zygmund} \|R_k g\|_{L^1} \leq C + C\int |g(x)|\ln^+|g(x)| \; \dd x. \end{equation}

By scaling, we can rewrite (\ref{zygmund}) as
\[ \|R_k g\|_{L^1} \leq \lambda + C\int |g(x)|\ln^+(|g(x)|/\lambda) \dd x \]
for any $\lambda>0$. For $|g|\geq\lambda$, we have $\ln^+(|g|/\lambda)=\ln|g|-\ln\lambda$. Then
\[ |\ln(|g|_{|_{|g|\geq\lambda}})| = \left| \ln\left(\frac{|g|_{|_{|g|\geq\lambda}}}{1+\|g\|_{L^q}}\right) + \ln(1+\|g\|_{L^q})\right| \leq \ln(1+\|g\|_{L^q}) + \left|\ln\left(\frac{|g|_{|_{|g|\geq\lambda}}}{1+\|g\|_{L^q}}\right)\right|. \]
Now assume $\lambda<1+\|g\|_{L^q}$ and take $x$ such that $\lambda\leq|g(x)|\leq \frac{(1+\|g\|_{L^q})^2}{\lambda}$. Then
\[ \left|\ln\left(\frac{|g(x)|}{1+\|g\|_{L^q}}\right) \right| \leq \left|\ln\left(\frac{\lambda}{1+\|g\|_{L^q}}\right)\right| \]
and in consequence
\[ |\ln|g(x)| \leq 2\ln(1+\|g\|_{L^q}) + |\ln\lambda|. \]
Choose $\lambda=\|g\|_{L^1}$. Then
\[\begin{aligned} \int_{\big\{|g|\leq \frac{(1+\|g\|_{L^q})^2}{\lambda}\big\}} |g|\ln^+(|g|/\lambda) \dd x &\leq \int_{\big\{ |g|\leq \frac{(1+\|g\|_{L^q})^2}{\lambda}\big\}} |g|(2\ln(1+\|g\|_{L^q})+|\ln\|g\|_{L^1}|) \dd x \\ 
&\leq \|g\|_{L^1}(2\ln(1+\|g\|_{L^q})+|\ln\|g\|_{L^1}|) \end{aligned}\]
What is left is the case $|g(x)|>\frac{(1+\|g\|_{L^q})^2}{\lambda}$. From the Chebyshev inequality, we have
\[ \left|\left\{x: |g(x)|>\frac{(1+\|g\|_{L^q})^2}{\|g\|_{L^1}} \right\}\right| \leq \left(\frac{\|g\|_{L^1}}{1+\|g\|_{L^q}}\right)^2. \]
Therefore from the H\"older inequality
\[ \int_{\big\{|g|>\frac{(1+\|g\|_{L^q})^2}{\lambda}\big\}} |g|\ln^+(|g|/\lambda)\dd x \leq \frac{\|g\|_{L^1}}{1+\|g\|_{L^q}}\left(\int |g|^2\ln(|g|/\lambda)^2 \dd x\right)^{1/2}. \]
Using the fact that both $\int |g|^2\ln|g| \dd x$ and $\int |g|^2\ln^2|g| \dd x$ are bounded by
\[ C\|g\|_{L^q}^{4-q}+\|g\|_{L^q}^q \]
for any $q>2$, we obtain
\[\begin{aligned} \int |g|^2\ln(|g|/\lambda)^2\dd x &= \int |g|^2\ln^2|g|\dd x - 2\ln\lambda\int |g|^2\ln|g| \dd x + |\ln\lambda|^2\int |g|^2 \dd x \\
&\leq (C\|g\|_{L^q}^{4-q}+\|g\|_{L^q}^q)(1+2|\ln\lambda|+|\ln\lambda|^2) \\
&= (C\|g\|_{L^q}^{4-q}+\|g\|_{L^q}^q)(1+|\ln\lambda|)^2. \end{aligned} \]
As $\frac{(Cs^{4-q}+s^q)^{1/2}}{1+s} \sim s^{\frac{q-2}{2}}$, after combining the estimates we get
\[ \int_{\big\{|g|>\frac{(1+\|g\|_{L^q})^2}{\lambda}\big\}} |g|\ln^+(|g|/\lambda)\dd x \leq C\|g\|_{L^1}(1+|\ln\|g\|_{L^1})\|g\|_{L^q}^{\frac{q-2}{2}}. \]
Putting all terms together,
\[ \|R_k g\|_{L^1} \leq C\|g\|_{L^1}\left(1 + \ln(1+\|g\|_{L^q})+|\ln\|g\|_{L^1}| + (1+|\ln\|g\|_{L^1}|)\|g\|_{L^q}^{\frac{q-2}{2}} \right) \]
and therefore we obtain inequality (\ref{mr_q}). \\
If $\int g \; \dd x \neq 0$, then we can apply this inequality to $\bar{g}(x)=g(x)-\frac{1}{|\mathbb{T}^d|}\int g \dd x$ and use the fact that the $L^p$ norms of $\bar{g}$ are bounded by norms of $g$ up to a constant. \qed

The same result holds if we replace $f$ by a composition of $f$ and the flow $x(t,y)$:
\begin{cor}\label{wniosek}
If $f$ and $g$ satisfy assumptions of Lemma \ref{log_ineq} and $x(t,y)$ is the regular Lagrangian flow of some $u(t,x)$ with bounded divergence, then
\[\begin{aligned} \left|\int f(x(t,y))g(y) \dd y \right| \leq & C\|f\|_{BMO}\|g\|_{L^1} \\
&\times \big(|\ln\|g\|_{L^1}| + \ln(e+\|g\|_{L^q}) + (1+|\ln\|g\|_{L^1}|)\|g\|_{L^q}^{\frac{q-2}{2}}\big). \end{aligned} \]
\end{cor}
Proof: We will first approximate $u$ with smooth vector fields, then make the change of variables and apply Lemma \ref{log_ineq}, and at the end show the convergence to the non-smooth case. \\
Let $J(t,y)$ be the Jacobian of $x$. By the properties of Lagrangian flows we have
\[ e^{-L} \leq J(t,y) \leq e^L, \text{ where } L=\int_0^T \|\textrm{div}u\|_\infty \dd t. \]
Now let us approximate $u$ by convolution, defining $u_\varepsilon$ as a convolution with standard convolution kernels in time and space. Then if $x_\varepsilon$ is a flow of $u_\varepsilon$, then $x_\varepsilon(t,\cdot)$ is the diffeomorphism and the Jacobian $J_\varepsilon$ of $x_\varepsilon$ still satisfies the bounds
\[ e^{-L} \leq J_\varepsilon(t,y) \leq e^L. \]
By the change of variables, we have
\[ \int f(x_\varepsilon(t,y))g(y)\dd y = \int f(x)\frac{g(y_\varepsilon(t,x))}{J_\varepsilon(t,y_\varepsilon(t,x))} \dd x, \]
where $y_\varepsilon(t,\cdot)=x_\varepsilon(t,\cdot)^{-1}$. Applying Lemma \ref{log_ineq}, we obtain inequality (\ref{mr_q}) but with $L^1$ and $L^q$ norms of $\frac{g(y_\varepsilon(t,\cdot))}{J_\varepsilon(t,y_\varepsilon(t,\cdot))}$ instead of $g$. However, changing the variables again and using the bounds on $J_\varepsilon$, we obtain for any $p\geq 1$
\[ \int \frac{|g(y_\varepsilon(t,x)|^p}{J_\varepsilon(t,y_\varepsilon(t,x))^p} \dd x = \int |g(y)|^p J_\varepsilon(t,y)^{1-p} \dd y \leq e^{(p-1)L}\int |g(y)|^p \dd y \]
and we are done. \\
Now we will show that indeed 
\begin{equation}\label{convergence} \int f(x_\varepsilon(t,y))g(y) \dd y \to \int f(x(t,y))g(y) \dd y \; \text{ with } \; \varepsilon\to 0. \end{equation}
By the stability of the flow, we have the pointwise convergence $x_\varepsilon(t,y)\to x(t,y)$ up to a subsequence. If $f\in C^\infty(\mathbb{T}^d)$, then (\ref{convergence}) holds by the dominated convergence theorem. Let us approximate $f$ by $f_\delta=f\ast\kappa_\delta$, where $\kappa_\delta$ is again the standard mollifier. As $f\in L^p$ for $p=q'$, we have
\[ \left|\int (f_\delta(x(t,y))-f(x(t,y))g(y) \dd y \right| \leq \|f_\delta(x(t,\cdot))-f(x(t,\cdot))\|_p\|g\|_q \]
and by the bounds on $J(t,y)$,
\[ \int |f_\delta(x(t,y))-f(x(t,y))|^p \; \dd x \leq e^{pL}\int |f_\delta(x)-f(x)|^p \; \dd x \to 0, \]
therefore we have the desired convergence. Moreover, by the Proposition \ref{prop_bmo} the norms $\|f_\delta\|_{BMO}$ in the right hand side of (\ref{mr_q}) are bounded by $\|f\|_{BMO}$, which ends the proof of the Corollary. \qed
 
\end{appendices}

 \bibliographystyle{abbrv}
\bibliography{biblio.bib}

\end{document}